\documentclass[12pt,a4paper]{article}
\usepackage{CJK,mathrsfs,amsfonts,geometry,fancyhdr,indentfirst}

\usepackage{amssymb}
\usepackage{mathrsfs}
\usepackage{amsthm}
\usepackage{amsmath}
\usepackage{graphicx}
\usepackage{color}
\usepackage{setspace}
\usepackage{titletoc}
\usepackage{tocloft}
\usepackage{appendix}

\newtheorem*{corollary*}{Corollary}

\theoremstyle{definition}

\newtheorem*{question*}{Question}

\begin{document}
\begin{CJK}{GBK}{song}
\begin{spacing}{1}
\newtheorem{remark}{remark}

\pagestyle{fancy}
\lhead{}
\chead{}
\rhead{}
\lfoot{}
\cfoot{}
\rfoot{\thepage}
\date{}

\title{Switch System}
\date{}
\author{Yuyi Zhang\footnote{Yuyi Zhang: 2021085001@fudan.edu.cn}}
\maketitle
\thispagestyle{empty}

\paragraph{Abstract:}
 Stability is a key property of dynamical systems. In some cases, we want to change unstable system into stable one to achieve certain goals in engineering. Here, we present an example of a $3$ dimensional switched system that alternates between unstable modes with the same period orbit, but exhibits stability for fast-switching frequencies and converges to the period orbit. Through mathematical analysis, we elucidate the stabilization mechanism and construct a general class of unstable systems that can form stable switched systems. The demonstration of this new system with period orbit provide a new sight in the study of switched system.
\paragraph{\textbf{[Key Words]}}\textbf{Switched System; Period Orbit; Stability}

\section{INTRODUCTION}
Stability is arguably the most fundamental property dynamical systems, which has been for decades. In switched systems, stability becomes priority due to engineering purposes. People adjust the switching frequencies when altering between stable systems in order to remain stability. In recent years, some examples about stable switched system generated from unstable mode have been made, opening a new door in this area.

However, researches on stability are largely focused on stable nodes, seldom other objects like period orbits have been discussed. In this paper, we present a stable fast-switching system that composes of two unstable system. More interestingly, the switched system converges to a period orbit, which is different to previous works. At the end, we provide a general type of this system, giving a way to generate stable switched systems with stable period orbits.

\section{Construction}
First, we construct two dynamic system with unstable period orbits. We will show that although period orbit of each system is unstable, it becomes stable if we switch between two system fast enough.

Consider the following two piecewise dynamic system, the first one is

\begin{equation}
\left\{
\begin{array}{lr}
\dot{x}=\frac{x}{\sqrt{x^2+y^2}}(-2z\sqrt{x^2+y^2}+10\sqrt{x^2+y^2})-y,&\\
\dot{y}=\frac{y}{\sqrt{x^2+y^2}}(-2z\sqrt{x^2+y^2}+10\sqrt{x^2+y^2})+x,& \sqrt{x^2+y^2}<\frac{1}{2}\\
\dot{z}=2z
\end{array}
\right.
\end{equation}

\begin{equation}
\left\{
\begin{array}{lr}
\dot{x}=\frac{x}{\sqrt{x^2+y^2}}(-10\sqrt{x^2+y^2}-z+10)-y,&\\
\dot{y}=\frac{y}{\sqrt{x^2+y^2}}(-10\sqrt{x^2+y^2}-z+10)+x,& \sqrt{x^2+y^2}\geq\frac{1}{2}\\
\dot{z}=2z
\end{array}
\right.
\end{equation}

while the second one is
\begin{equation}
\left\{
\begin{array}{lr}
\dot{x}=\frac{x}{\sqrt{x^2+y^2}}(-2\sqrt{x^2+y^2}+2z\sqrt{x^2+y^2})-y,&\\
\dot{y}=\frac{y}{\sqrt{x^2+y^2}}(-2\sqrt{x^2+y^2}+2z\sqrt{x^2+y^2})+x,& \sqrt{x^2+y^2}<\frac{1}{2}\\
\dot{z}=-10z
\end{array}
\right.
\end{equation}

\begin{equation}
\left\{
\begin{array}{lr}
\dot{x}=\frac{x}{\sqrt{x^2+y^2}}(2\sqrt{x^2+y^2}+z-2)-y,&\\
\dot{y}=\frac{y}{\sqrt{x^2+y^2}}(2\sqrt{x^2+y^2}+z-2)+x,& \sqrt{x^2+y^2}\geq\frac{1}{2}\\
\dot{z}=-10z
\end{array}
\right.
\end{equation}

Clearly, each of them forms a continuous vector field over the whole space. In other words, each of them is indeed a continuous $3$-dimension dynamic system. Also, they have the same period orbit $x^2+y^2=1, z=0$, which can be checked easily.

A more detailed analysis shows that period orbit of the first system is unstable, as the derivative at $z$-direction never converges to zero if $z\neq 0$.

\begin{figure}[htbp]
\centering
\includegraphics[width=0.6\textwidth]{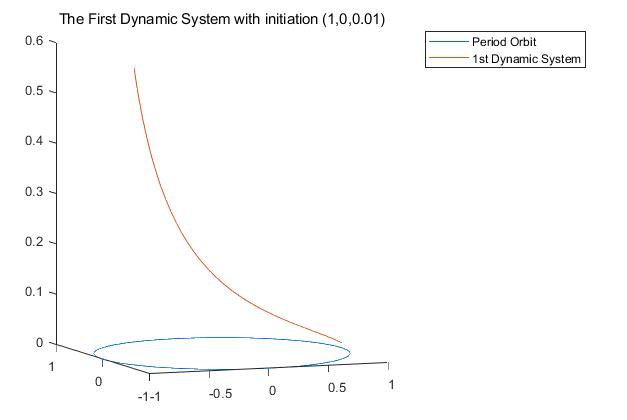}\\
\caption{The blue line represents period orbit $x^2+y^2=1,z=0$; the red line represents trajectory of the 1st dynamic system with given initiation.}
\end{figure}

Meanwhile, we can show the similar result for the second system if we analyse on a two-dimension section of the origin system. We may choose $XOZ$-plane for simplicity, then the system reduces to a two-dimension system.
\begin{equation}
\left\{
\begin{array}{lr}
\dot{x}=-2x+2xz,& 0\geq x<\frac{1}{2}\\
\dot{z}=-10z,
\end{array}
\right.
\end{equation}

\begin{equation}
\left[
\begin{matrix}
\dot{x} \\
\dot{z}
\end{matrix}
\right]
=
\left[
\begin{matrix}
     2 & 1 \\
     0 & -10
\end{matrix}
\right]
\left[
\begin{matrix}
     x-1 \\
     z
\end{matrix}
\right],x\geq \frac{1}{2}
\end{equation}

A standard analysis on two-dimension system shows that it's unstable since a eigenvalue of matrix in equation $(6)$ has positive real part. We can also draw out its trajectory as follow.
\begin{figure}[htbp]
\centering
\includegraphics[width=0.6\textwidth]{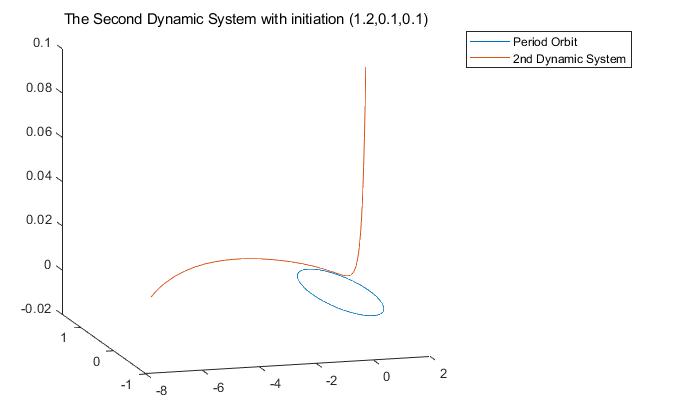}\\
\caption{The blue line represents period orbit $x^2+y^2=1,z=0$; the red line represents trajectory of the 2nd dynamic system with given initiation.}
\end{figure}

\section{Switching System}

Although both system are unstable, they form a stable period orbit under a fast switching system.
\begin{figure}[htbp]
\centering
\includegraphics[width=0.6\textwidth]{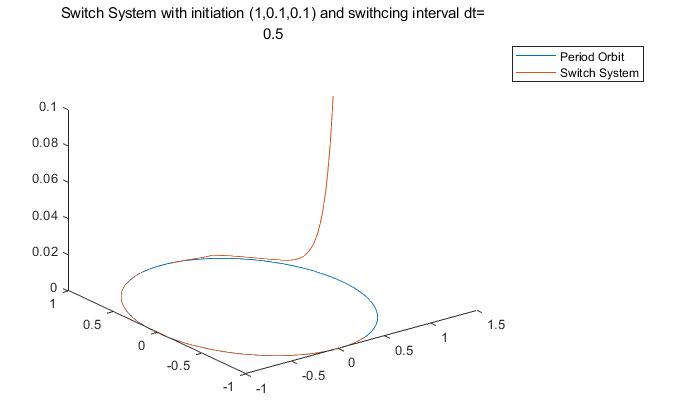}\\
\caption{The blue line represents period orbit $x^2+y^2=1,z=0$; the red line represents trajectory of the switching dynamic system with given initiation; we switch from one system to another for every $dt=0.5$ second.}
\end{figure}

Clearly, the trajectory of switching system converges to the period orbit. This can be proved by average system, as their average is the following.

\begin{equation}
\left\{
\begin{array}{lr}
\dot{x}=\frac{x}{\sqrt{x^2+y^2}}(4\sqrt{x^2+y^2})-y,&\\
\dot{y}=\frac{y}{\sqrt{x^2+y^2}}(4\sqrt{x^2+y^2})+x,& \sqrt{x^2+y^2}<\frac{1}{2}\\
\dot{z}=-4z
\end{array}
\right.
\end{equation}

\begin{equation}
\left\{
\begin{array}{lr}
\dot{x}=\frac{x}{\sqrt{x^2+y^2}}(-4\sqrt{x^2+y^2}+4)-y,&\\
\dot{y}=\frac{y}{\sqrt{x^2+y^2}}(-4\sqrt{x^2+y^2}+4)+x,& \sqrt{x^2+y^2}\geq\frac{1}{2}\\
\dot{z}=-4z
\end{array}
\right.
\end{equation}

Using the same trick, we reduce the system to the section $XOZ$-plane. The reducing system is
\begin{equation}
\left\{
\begin{array}{lr}
\dot{x}=4x,& 0\geq x<\frac{1}{2}\\
\dot{z}=-4z,
\end{array}
\right.
\end{equation}

\begin{equation}
\left[
\begin{matrix}
\dot{x} \\
\dot{z}
\end{matrix}
\right]
=
\left[
\begin{matrix}
     -4 & 0 \\
     0 & -4
\end{matrix}
\right]
\left[
\begin{matrix}
     x-1 \\
     z
\end{matrix}
\right],x\geq \frac{1}{2}.
\end{equation}
Since the matrix has only negative eigenvalue, the period orbit is stable. Note that stability of period orbit is a local property, which means we only need to consider points that are close enough to the orbit. Hence, it's enough to analyse the second part($>\frac{1}{2}$) since the orbit lies in the second part for every system. Therefore our consideration is valid.

It's worth noting that switching system remains unstable if switching pace is too slow as we can see from the following.
\begin{figure}[htbp]
\centering
\includegraphics[width=0.6\textwidth]{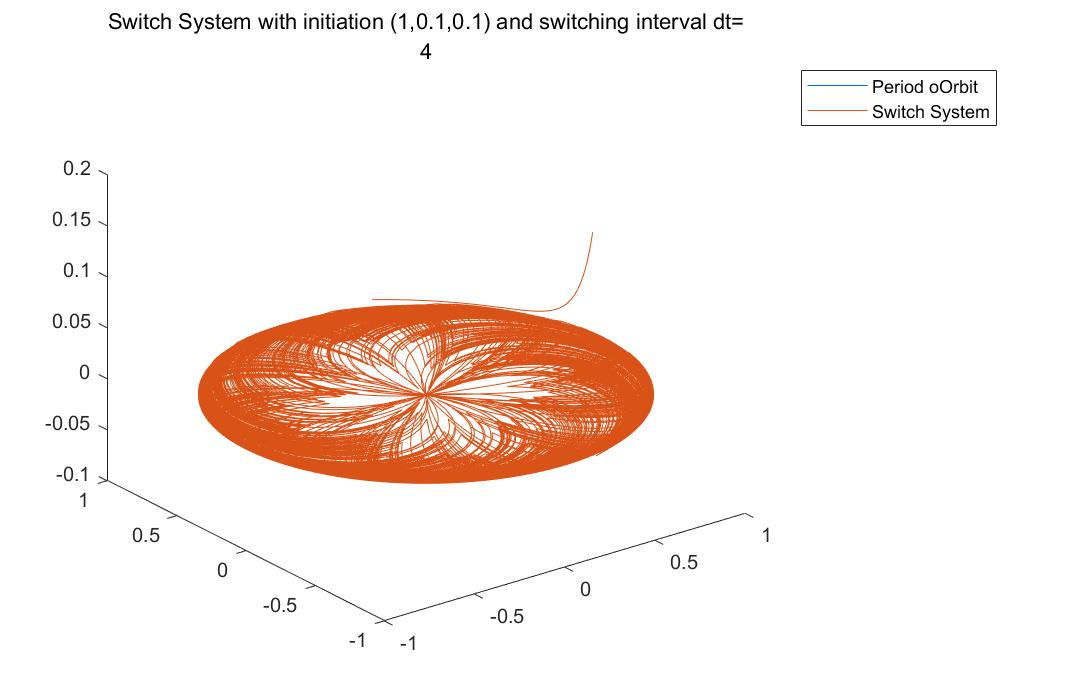}\\
\caption{The blue line represents period orbit $x^2+y^2=1,z=0$; the red line represents trajectory of the switching dynamic system with given initiation; we switch from one system to another for every $dt=4$ second. Clearly, although the trajectory performs certain pattern, it never converges to the period orbit, which means the orbit is unstable in switching system.}
\end{figure}

\section{Mathematical Analysis with Cylindrical Coordinates}
Under cylindrical coordinates, all systems above become much more clear. Let $x=rcos\theta, y=rsin\theta$, system $(1)-(2)$ turns into
\begin{equation}
\left\{
\begin{array}{lr}
\dot{r}=10r-2rz,\\
\dot{\theta}=1,& 0\geq r<\frac{1}{2} \\
\dot{z}=2z,
\end{array}
\right.
\end{equation}

\begin{equation}
\left[
\begin{matrix}
\dot{r} \\
\dot{\theta}\\
\dot{z}
\end{matrix}
\right]
=
\left[
\begin{matrix}
     -10 & 0 &-1 \\
     0 & 0 & 0 \\
     0 & 0 & 2
\end{matrix}
\right]
\left[
\begin{matrix}
     r-1 \\
     \theta\\
     z
\end{matrix}
\right]
+
\left[
\begin{matrix}
    0 \\
    1 \\
    0
\end{matrix}
\right],r\geq \frac{1}{2}.
\end{equation}

and system $(3)-(4)$ turns into
\begin{equation}
\left\{
\begin{array}{lr}
\dot{r}=-2r+2rz,\\
\dot{\theta}=1,& 0\geq r<\frac{1}{2} \\
\dot{z}=-10z,
\end{array}
\right.
\end{equation}

\begin{equation}
\left[
\begin{matrix}
\dot{r} \\
\dot{\theta}\\
\dot{z}
\end{matrix}
\right]
=
\left[
\begin{matrix}
     2 & 0 &1 \\
     0 & 0 & 0 \\
     0 & 0 & -10
\end{matrix}
\right]
\left[
\begin{matrix}
     r-1 \\
     \theta\\
     z
\end{matrix}
\right]
+
\left[
\begin{matrix}
    0 \\
    1 \\
    0
\end{matrix}
\right],r\geq \frac{1}{2}.
\end{equation}

By the same reason above, we only need to consider the second part for each system. And each of them is a linear system with a eigenvalue that has positive real part, thus their period orbits $r=1, z=0$ can't be stable, which is presented by stimulation in section $1$.

Also, their average system becomes
\begin{equation}
\left\{
\begin{array}{lr}
\dot{r}=4r,\\
\dot{\theta}=1,& 0\geq r<\frac{1}{2} \\
\dot{z}=-4z,
\end{array}
\right.
\end{equation}

\begin{equation}
\left[
\begin{matrix}
\dot{r} \\
\dot{\theta}\\
\dot{z}
\end{matrix}
\right]
=
\left[
\begin{matrix}
     -4 & 0 &0 \\
     0 & 0 & 0 \\
     0 & 0 & -4
\end{matrix}
\right]
\left[
\begin{matrix}
     r-1 \\
     \theta\\
     z
\end{matrix}
\right]
+
\left[
\begin{matrix}
    0 \\
    1 \\
    0
\end{matrix}
\right],r\geq \frac{1}{2}.
\end{equation}
Clear, the second part of equation is a linear system with only negative real eigenvalue, hence its period orbit is stable.

In general, we can construct a class of system with unstable period orbits, which can form fast-switching systems under some mild conditions. For simplicity, we restrict our attention in $3$ dimensional real case. Consider systems $\{S_i\}_{I}$ with following forms.
\begin{equation}
\left\{
\begin{array}{lr}
\dot{r}=-a_ir+2zb_ir,\\
\dot{\theta}=1,& 0\geq r<\frac{1}{2}d \\
\dot{z}=c_iz,
\end{array}
\right.
\end{equation}

\begin{equation}
\left[
\begin{matrix}
\dot{r} \\
\dot{\theta}\\
\dot{z}
\end{matrix}
\right]
=
\left[
\begin{matrix}
     a_i & 0 &b_i \\
     0 & 0 & 0 \\
     0 & 0 & c_i
\end{matrix}
\right]
\left[
\begin{matrix}
     r-d \\
     \theta\\
     z
\end{matrix}
\right]
+
\left[
\begin{matrix}
    0 \\
    1 \\
    0
\end{matrix}
\right],r\geq \frac{1}{2}d.
\end{equation}

Obviously, it has period orbit$r=d,z=0$. If $\sum_I a_i<-1,\sum_I c_i<-1$ and $\sum_I b_i=0$, then the average system of this class of systems has stable period orbit(condition on $b_i$ is not necessary since it doesn't affect eigenvalues). We may choose some of $a_i,c_i$ to be positive to make each system unstable, and our construction is complete.

\section{conclusion}
Stabilization of unstable modes via periodic or stochastic perturbations is a critical problem in the study of dynamical systems, with applications in physical, biological, and engineering systems. To date, there was no example of stabilization of a switched system with stable period orbit composed of only unstable systems through periodic or stochastic switching. Here we provide this new phenomenon, making switched system more fascinating. The key point is that although each subsystems may be unstable, their average system is stable, forcing the switched system trapped in a bounded region.

The simplicity and richness of our system offers a good example finding more complicated stable dynamics, not just having stable or unstable nodes, but something like period orbits and even more.

\end{spacing}

\end{CJK}
\end{document}